\documentclass[12pt]{amsart}
\usepackage{amssymb}

\newtheorem{theorem}{Theorem}

\newcommand{\rf}[1]{{\rm (}\ref{#1}{\rm )}}
\def\qed{{\hfill$\Box$}}
\def\RN{\mathbb{R}^N}

\title[Nonlinear equation with fractional Laplacian]{\bf  Decay of mass for
 nonlinear  equation \\  with fractional Laplacian}

\author{Ahmad Fino}
\address{
La\-bo\-ra\-toire MIA et 
D\'e\-par\-te\-ment de Math\'ematiques,
Universit\'e de La Rochelle, Avenue Michel Cr\'epeau,
17042 La Rochelle Cedex, France
\&
LaMA-Liban, Lebanese University,
P.O. Box 826, Tripoli, Lebanon 
}
\email{afino01@univ-lr.fr}

\author{Grzegorz Karch}
\address{Instytut Matematyczny, Uniwersytet Wroc\l awski,
 pl.~Grun\-wal\-dzki 2/4, 50-384 Wroc\l aw, Poland}
\email{karch{@}math.uni.wroc.pl}
\urladdr{http://www.math.uni.wroc.pl/$\sim$karch}

\date{\today}



\thanks{
The preparation of this paper was supported in part by the European
Commission Marie Curie Host Fellowship for the Transfer of Knowledge
``Harmonic Analysis, Nonlinear Analysis and Probability''
MTKD-CT-2004-013389. The first author gratefully thanks the
Mathematical Institute of Wroc\l aw University for the warm
hospitality. The preparation of this paper by the second author was
also partially supported by the MNiSW grant  N201 022 32 / 09 02. }

\begin{document}
\baselineskip=16pt

\begin{abstract}
The large time
behavior of nonnegative solutions to the reaction-diffusion equation
$\partial_t u=-(-\Delta)^{\alpha/2}u - u^p,$
$(\alpha\in(0,2], \;p>1)$  posed on
$\mathbb{R}^N$ and supplemented with an integrable initial condition is
studied. We show that the
anomalous diffusion term  determines the large time asymptotics for
$p>1+{\alpha}/{N},$ while nonlinear effects win if $p\leq1+{\alpha}/{N}.$
\end{abstract}

\maketitle

\section{Introduction}

We study the behavior, as $t\to\infty$, of solutions to  the
following initial value problem for the  reaction-diffusion equation with
the anomalous diffusion
\begin{eqnarray}\label{eq}
\partial_t u &=& -\Lambda^\alpha u  + \lambda u^p,\qquad x\in\mathbb{R}^N,t>0,\\
 u(x,0) &=& u_0(x),\label{ini}
\end{eqnarray}
where the  pseudo-differential operator
$\Lambda^\alpha=(-\Delta)^{\alpha/2}$ with $0<\alpha\leq2$  is
 defined by the Fourier transformation:
$
\widehat{\Lambda^\alpha u}(\xi)=|\xi|^\alpha\widehat{u}(\xi).
$
Moreover, we assume that  $\lambda\in\{-1,1\}$ and $p>1$.

Nonlinear evolution problems involving fractional Laplacian
describing {\it  the ano\-ma\-lous diffusion} (or $\alpha$-stable
L\'evy diffusion) have been extensively studied in the mathematical
and physical literature (see \cite{BKW01, KW08, 6} for references).
One of possible ways to understand the interaction between  the
anomalous diffusion operator (given by $\Lambda^\alpha$ or, more
generally, by the L\'evy diffusion operator) and the nonlinearity in
the equation under consideration is the study of the large time
asymptotics of solutions   to such equations. Our goal is  to
contribute to this theory and our results can be summarized as
follows. For $\lambda=-1$ in equation \rf{eq}, nonnegative solutions
to the Cauchy problem exist globally in time.  Hence, we  study the
decay properties of the mass $M(t)= \int_{\mathbb{R}^N}u(x,t)\,dx$
of the solutions $u=u(x,t)$ to problem \rf{eq}-\rf{ini}. We prove
that $ \lim_{t\to\infty}M(t)=M_\infty >0$ for $p>1+{\alpha}/{N}$
(cf. Theorem 1, below), while
 $M(t)$ tends to zero as $t\to\infty$
if $1<p\leq 1+\alpha/N$ (cf. Theorem 2).
As a by-product of our analysis, we show the blow-up of all
nonnegative solutions to \rf{eq}-\rf{ini} with $\lambda=1$ in the case of the
critical nonlinearity exponent $p=1+\alpha/N$ (see Theorem 3, below).

The idea which allows to express the competition between  diffusive and  nonlinear terms
in an evolution equation by studying the large time behavior of the space integral
of a solution was already introduced by Ben-Artzi \& Koch \cite{BK99}
who considered
the viscous Hamilton-Jacobi equation $u_t=\Delta u-|\nabla u|^p$ (see also Pinsky \cite{21}).
An analogous result
for the equation $u_t=\Delta u+|\nabla u|^p$
(with the growing-in-time mass of solutions) was proved
by Lauren\c cot \& Souplet \cite{LS03}. Such questions concerning
the asymptotic behavior of solutions to the  Hamilton-Jacobi equation
with the L\'evy diffusion operator were answered in \cite{KW08}.

In the case of the classical reaction-diffusion equation
({\it i.e.}~equation \rf{eq} with $\alpha=2$), for $p<1+{2}/{N}$,
Fujita \cite{8} proved the nonexistence of nonnegative global-in-time
solution for any nontrivial initial condition.
On other hand,
if $p>1+{2}/{N},$ global solutions do exist for any
sufficiently small nonnegative initial data.
The proof of a blow-up of all nonnegative solutions in the critical case
$p=1+{2}/{N}$ was completed in \cite{11,22,17}.
Analogous blow-up results for problem \rf{eq}-\rf{ini}
 with the fractional Laplacian
(and with the critical exponent $p=1+\alpha/N$
for the existence/nonexistence of solutions)
are contained {\it e.g.}  in \cite{22,9,10, BLW02}.

\section{Statement of results}

In all theorems below, we always assume that $u=u(x,t)$ is the nonnegative
(possibly weak)
solution of problem \rf{eq}-\rf{ini} corresponding to the nonnegative initial datum
$u_0\in L^1(\mathbb{R}^N).$ Let   $u_0\not\equiv 0$, for simplicity of the exposition.
We refer the reader to \cite{6} for several results on the existence, the uniqueness
and the regularity of solutions  to \rf{eq}-\rf{ini} as well as for the proof of the maximum principle
(which assures that the solution is nonnegative if the corresponding
initial datum is so).

First, we deal with the equation \rf{eq} containing the absorbing
nonlinearity $(\lambda=-1)$ and we study the decay of the ``mass''
\begin{equation}\label{mass}
    M(t)\equiv \int_{\mathbb{R}^N}u(x,t)\,dx=\int_{\mathbb{R}^N}u_0(x)\,dx -\int_0^t\int_{\mathbb{R}^N}u^p(x,s)\,dxds.
\end{equation}

\noindent{\bf Remark.~}In order to obtain equality \rf{mass}, it
suffices to integrate equation \rf{eq} with respect to $x$ and $t$.
Another method which leads to \rf{mass} and which requires weaker
regularity assumptions on a solution consists in integrating with
respect to $x$ the integral formulation of problem \rf{eq}-\rf{ini}
(see \rf{3.8}, below) and using the Fubini theorem. \qed

\bigskip

Since we limit ourselves to nonnegative solutions, the function
$M(t)$ defined in \rf{mass} is nonnegative and non-increasing. Hence, the limit
$M_\infty = \lim_{t\rightarrow\infty}M(t)$ exists and we answer
the question whether it is equal to zero or not.

In our first theorem, the  diffusion phenomena determine the large time
asymptotics of solutions to \rf{eq}-\rf{ini}.

\begin{theorem}
Assume that $u=u(x,t)$ is a nonnegative nontrivial solution of
\rf{eq}-\rf{ini} with $\lambda=-1$ and $p>1+{\alpha}/{N}.$ Then
$
\lim_{t\to\infty}M(t)=M_\infty >0.
$

Moreover, for all $q\in[1,\infty)$
\begin{equation}\label{2.3}
    t^{\frac{N}{\alpha}\left(1-\frac{1}{q}\right)}\|u(t)-M_\infty
P_\alpha(t)\|_q\to
0\quad\hbox{as} t\to\infty,
\end{equation}
where  the function $P_\alpha(x,t)$ denotes the fundamental
solution of the linear equation $u_t + \Lambda^\alpha u = 0$
(cf. equation \rf{3.3} below).
\end{theorem}

In the remaining range of $p$, the mass $M(t)$ converges to zero and
this phenomena can be interpreted as the domination of nonlinear
effects in the large time asymptotic of solutions to \rf{eq}-\rf{ini}.
Note here that the mass $M(t)=\int_{\RN}u(x,t)\,dx$
 of every solution to linear equation
$u_t+\Lambda^\alpha u=0$ is constant in time.

\begin{theorem}
Assume that $u=u(x,t)$ is a nonnegative solution of problem
\rf{eq}-\rf{ini} with $\lambda=-1$ and $1< p\leq1+{\alpha}/{N}.$
Then
$
\lim_{t\rightarrow\infty}M(t)=0.
$
\end{theorem}

Let us emphasize that the proof of Theorem 2 is based on the
so-called the rescaled test function method which was
 used  by Mitidieri \& Pokhozhaev ({\it cf.~e.g.}~\cite{19,20} and
the references therein)   to prove
 the nonexistence of solutions to  nonlinear elliptic and parabolic equations.

As the by-product of our analysis, we can also contribute to the
theory on the blow-up of solutions to \rf{eq}-\rf{ini} with
$\lambda=+1.$ Recall that the method of the rescaled test function
(which we also apply here)  was use in
\cite{9,10} to show the blow-up of all positive
solutions to \rf{eq}-\rf{ini}  with $\lambda=1$ and
$p<1+{\alpha}/{N}.$ Here, we complete
that result by the simple proof of the blow-up in
the critical case $p=1+{\alpha}/{N}.$

\begin{theorem}
If  $\lambda=1,$ $\alpha \in (0,2]$ and $ p = 1+{\alpha}/{N},$
then any nonnegative nonzero solution of \rf{eq}-\rf{ini}
 blows up
in a finite time.
\end{theorem}


\section{Proofs of Theorems 1, 2, and 3}

Note first that any (sufficiently regular) nonnegative solution to \rf{eq}-\rf{ini}
 satisfies
\begin{equation}\label{3.1}
    0\leq\int_{\mathbb{R}^N}u(x,t)\,dx=\int_{\mathbb{R}^N}u_0(x)\,dx + \lambda\int_0^t\int_{\mathbb{R}^N}u^p(x,s)\,dx\,ds.
\end{equation}
\noindent Hence, for $\lambda=-1$ and $u_0\in L^1(\mathbb{R}^N),$ we
immediately obtain
\begin{equation}\label{3.2}
u\in L^\infty([0,\infty),L^1(\mathbb{R}^N))\cap
L^p(\mathbb{R}^N\times(0,\infty)).
\end{equation}

\bigskip
{\it Proof of Theorem 1.}
First, we recall that the fundamental solution
$P_\alpha=P_\alpha(x,t)$ of the linear equation $\partial_t
u+\Lambda^\alpha u=0$ can be written  via the Fourier transform
as follows
\begin{equation}\label{3.3}
P_\alpha(x,t)=t^{-N/\alpha}P_\alpha(xt^{-1/\alpha},1)=\frac{1}{(2\pi)^{N/2}}\int_{\mathbb{R}^N}e^{ix.\xi-t|\xi|^\alpha}\,d\xi.
\end{equation}
It is well-known that for each $\alpha\in(0,2],$ this function
satisfies
\begin{equation}\label{3.4}
    P_\alpha(1)\in L^\infty(\mathbb{R}^N)\cap
L^1(\mathbb{R}^N),\quad
P_\alpha(x,t)\geq0,\quad\int_{\mathbb{R}^N}P_\alpha(x,t)\,dx=1,
\end{equation}
\noindent for all $x\in\mathbb{R}^N$ and $t>0.$ Hence, using the
Young inequality for the convolution and the self-similar form of
$P_\alpha,$ we have
\begin{eqnarray}
\|P_\alpha(t)\ast u_0\|_p&\leq&Ct^{-N(1-1/p)/\alpha}\|u_0\|_1,\label{3.5}\\
\|\nabla P_\alpha(t)\|_p&=&Ct^{-N(1-1/p)/\alpha-1/\alpha},\label{3.6}\\
\|P_\alpha(t)\ast u_0\|_p&\leq& \|u_0\|_p,\label{3.7}
\end{eqnarray}
for all $p\in[1,\infty]$ and $t>0.$

In the next step, using the following well-known integral  representation of
solutions to \rf{eq}-\rf{ini}
\begin{equation}\label{3.8}
    u(t)= P_\alpha(t)\ast u_0 -\int_0^tP_\alpha(t-s)\ast
u^p(s)\,ds,
\end{equation}
 we immediately obtain the estimate $0\leq u(x,t)\leq
P_\alpha(x,t)\ast u_0(x).$ Hence, by \rf{3.5} and \rf{3.7} we get
\begin{eqnarray}\label{3.9}
  \|u(t)\|^p_p &\leq& \|P_\alpha(t)\ast u_0\|_p^p\nonumber\\
  &\leq& \min\left\{C t^{-N(p-1)/\alpha}\|u_0\|_1^p ; \|u_0\|_p^p\right\}\equiv H(t,p,\alpha,u_0).
\end{eqnarray}

Now, for fixed $\varepsilon\in(0,1],$ we consider the solution
$u^\varepsilon=u^\varepsilon(x,t)$ of \rf{eq}-\rf{ini}  with the initial
condition $\varepsilon u_0(x).$ The comparison principle implies
that $0\leq u^\varepsilon(x,t)\leq u(x,t)$ for every
$x\in\mathbb{R}^N$ and $t>0.$ Hence, it suffices to show that for
small $\varepsilon>0,$ which will be determined later, we have
\begin{equation*}\label{3.10}
M^\varepsilon_\infty\equiv\lim_{t\rightarrow\infty}\int_{\mathbb{R}^N}u^\varepsilon(x,t)\,dx >0.
\end{equation*}

Note first the using  equality $(\ref{3.1})$ in the case of the
solution $u^\varepsilon,$ we obtain
\begin{equation}\label{3.11}
  M^\varepsilon_\infty = \varepsilon\left\{\int_{\mathbb{R}^N}u_0(x)\,dx -\frac{1}{\varepsilon} \int_0^\infty\int_{\mathbb{R}^N}\left(u^\varepsilon(x,t)\right)^p\,dx\,dt\right\}.
\end{equation}
Now, we apply $(\ref{3.9})$ with $u$ replaced by
$u^\varepsilon$. Observe that the function $H$ defined in
\rf{3.9} satisfies
$H(t,p,\alpha,\varepsilon
u_0) = \varepsilon^p H(t,p,\alpha,u_0).$ Hence
\begin{eqnarray*}
  \frac{1}{\varepsilon}\int_0^\infty\int_{\mathbb{R}^N}
\left(u^\varepsilon (x,t)\right)^p\,dxdt
&\leq&\frac{1}{\varepsilon}\int_0^\infty H(t,p,\alpha,\varepsilon u_0)\,dt \\
  &=&\varepsilon^{1-p}\int_0^\infty H(t,p,\alpha,u_0)\,dt.
\end{eqnarray*}
 It is follows immediately from the definition of the
function $H$
that the integral on the right-hand side is convergent
for $p>1+{\alpha}/{N}.$ Consequently,
$$
\frac{1}{\varepsilon} \int_0^\infty\int_{\mathbb{R}^N}\left(u^\varepsilon(x,t)\right)^p\,dx\,dt\to
0\qquad\hbox{as}\quad\varepsilon\searrow 0,
$$
and the constant $M^\varepsilon_\infty$ given by
$(\ref{3.11})$ is positive for sufficiently small $\varepsilon>0.$

From now on, the proof of the asymptotic relation \rf{2.3} is standard,
hence, we shall be brief in details. First we recall that for every
$u_0\in L^1(\mathbb{R}^N)$ we have
\begin{equation}\label{L1:lin}
\lim_{t\to\infty} \|P_\alpha(t)*u_0-MP_\alpha(t)\|_1=0,
\end{equation}
where $M=\int_{\mathbb{R}^N}u_0(x)\,dx$. This is the immediate consequence
of the Taylor argument combined with an
approximation argument. Details of this reasoning can be found  in
\cite[Lemma 3.3]{BKW01}.

Now, to complete the proof of Theorem 1, we adopt the reasoning from \cite{LS03}.
 It follows from the
integral equation $(\ref{3.8})$ and inequality $(3.7)$ with $p=1$ that
\begin{equation*}\label{3.14}
\|u(t)-P_\alpha(t-t_0)\ast u(t_0)\|_1\leq \int_{t_0}^t
\|u(s)\|_p^p\,ds\quad\hbox{for all}\quad t\geq t_0\geq 0.
\end{equation*}
Hence, using the triangle inequality we infer
\begin{equation}\label{3.15}
\begin{split}
  \|u(t)-M_\infty P_\alpha(t)\|_1 \leq& \int_{t_0}^t
\|u(s)\|_p^p\,ds\\
&+\|P_\alpha(t-t_0)\ast
u(t_0)-M(t_0)P_\alpha(t-t_0)\|_1\\
&+\|M(t_0)(P_\alpha(t-t_0)-P_\alpha(t))\|_1\\
&+
\|P_\alpha(t)\|_1\left|M(t_0)-M_\infty\right|.
\end{split}
\end{equation}
Applying  first \rf{3.4} and \rf{L1:lin} with $u_0=u(t_0)$, and next
passing to the limit as $t\to\infty$ on the right-hand side of
\rf{3.15}, we obtain
$$\limsup_{t\to\infty}\|u(t) - M_\infty
P_\alpha(t)\|_1\leq\int_{t_0}^\infty
\|u(s)\|_p^p\,ds + \left|M(t_0)-M_\infty\right|.$$\\
 By letting $t_0$ go to $+\infty$ and using \rf{3.2} we conclude
 that
\begin{equation}\label{3.16}
\|u(t) - M_\infty
P_\alpha(t)\|_1\to 0\quad\quad\hbox{as} \quad
t\to\infty.
\end{equation}

In order to obtain the asymptotic term for $p>1$, observe that by the
integral equation \rf{3.8} and estimate  \rf{3.5},
 for each $m\in[1,\infty],$ we have
\begin{equation}\label{3.17}
    \|u(t)\|_m\leq\|P_\alpha(t)\ast
    u_0\|_m\leq Ct^{-N(1-1/m)/\alpha}\|u_0\|_1.
\end{equation}
 Hence, for every $q\in[1,m),$ using the H\"older inequality,
we  obtain
\begin{equation*}\label{3.18}
\begin{split}
\|u(t)-M_\infty P_\alpha(t)\|_q&\leq \|u(t)-M_\infty
P_\alpha(t)\|_1^{1-\delta}\left(\|u(t)\|_m^\delta + \|M_\infty P_\alpha(t)\|_m^\delta\right)\\
&\leq Ct^{-N(1-1/q)/\alpha}\|u(t)-M_\infty
P_\alpha(t)\|_1^{1-\delta},
\end{split}
\end{equation*}
 with $\delta=(1-1/q)/(1-1/m).$
Finally, applying \rf{3.16} we complete the proof of Theorem~1.
\qed


\bigskip

{\it Proof of Theorem 2.}
Let us define  the  function
$\varphi(x,t)=\left(\varphi_1(x)\right)^\ell\left(\varphi_2(t)\right)^\ell$
where
$$
\ell=\frac{2p-1}{p-1},\quad
\varphi_1(x)=\psi\left(\frac{|x|}{BR}\right), \quad
\varphi_2(t)=\psi\left(\frac{t}{R^\alpha}\right), \quad R>0,
$$
 and
$\psi$ is a smooth non-increasing function on $[0,\infty)$ such that
$$\psi(r)=\left\{\begin {array}{l}1\qquad \quad \mbox {if }0\leq r\leq 1,\\
0\qquad \quad \mbox {if }r\geq 2.
\end {array}\right.$$
The constant $B>0$ in the definition of $\varphi_1$ is  fixed and
will be chosen later. In
fact, it plays some role in the critical case $p=1+{\alpha}/{N}$
only while in the subcritical case $p<1+\alpha/N$ we simply put
$B=1$.
 In the following, we denote by $\Omega_1$ and $\Omega_2$
the supports of $\varphi_1$ and $\varphi_2,$ respectively:
\begin{equation*}\label{3.19}
\Omega_1=\left\{x\in\mathbb{R}^N \;:\; |x|\leq
2BR\right\},\qquad\Omega_2=\left\{t\in[0,\infty)\;:\;  t\leq
2R^\alpha\right\}.
\end{equation*}

Now, we multiply equation \rf{eq}  by
$\varphi(x,t)$ and integrate with respect to $x$ and $t$ to obtain
\begin{eqnarray}
&&\hspace{-2cm}\int_{\Omega_1}u_0(x)\varphi(x,0)\,dx -\int_{\Omega_2}\int_{\Omega_1}u^p(x,t) \varphi(x,t)\,dxdt\nonumber\\
&=&\int_{\Omega_2}\int_{\mathbb{R}^N}u(x,t)\varphi_2(t)^\ell \Lambda^\alpha(\varphi_1(x))^\ell\,dxdt\nonumber\\
&&-\int_{\Omega_2}\int_{\Omega_1}u(x,t) \varphi_1(x)^\ell\partial_t\varphi_2(t)^\ell\,dxdt\label{3.20}\\
&\leq& \ell\int_{\Omega_2}\int_{\Omega_1}u(x,t)\varphi_2(t)^\ell\varphi_1(x)^{\ell-1}\Lambda^\alpha\varphi_1(x)\,dxdt\nonumber\\
&&- \ell\int_{\Omega_2}\int_{\Omega_1}u(x,t) \varphi_1(x)^\ell\varphi_2(t)^{\ell-1}\partial_t\varphi_2(t)\,dxdt.\nonumber
\end{eqnarray}
In \rf{3.20}, we have used the inequality $ \Lambda^\alpha
\varphi_1^\ell \leq\ell\varphi_1^{\ell-1}\Lambda^\alpha\varphi_1, $
(see \cite[Prop. 2.3]{5} and \cite[Prop. 3.3]{13} for its proof)
which is valid for all $\alpha\in(0,2]$,  $\ell\geq 1,$ and any
sufficiently regular, nonnegative, decaying at infinity function
$\varphi_1$.

Hence, by the $\varepsilon$-Young inequality
$ab\leq \varepsilon
a^p + C(\varepsilon)b^{\ell -1}$ (note that $1/p+1/(\ell-1)=1$)
with   $\varepsilon>0,$ we deduce from \rf{3.20}
\begin{equation}\label{3.21}
\begin{split}
&\int_{\Omega_1}u_0(x)\varphi(x,0)\,dx-(1+2\ell\varepsilon)
\int_{\Omega_2}\int_{\Omega_1}u^p(x,t) \varphi(x,t)\,dx\,dt\\
 & \leq
 C(\varepsilon)\ell \left\{\int_{\Omega_2}\int_{\Omega_1}\varphi_1
\varphi_2^{\ell}\left|\Lambda^\alpha \varphi_1\right|^{\ell-1}\,dxdt
+\int_{\Omega_2}\int_{\Omega_1}\varphi_1^{\ell}\varphi_2
\left|\partial_t\varphi_2\right|^{\ell -1}\,dxdt\right\}.
\end{split}
\end{equation}

 Recall now that the functions $\varphi_1$ and $\varphi_2$
depend on $R>0.$ Hence changing the variables $\xi=R^{-1}x$ and
$\tau=R^{-\alpha}t,$ we easily obtain from \rf{3.21} the
following estimate
\begin{equation}\label{3.22}
    \int_{\Omega_1}u_0(x)\varphi(x,0)\,dx -(1+2\ell \varepsilon)\int_{\Omega_2}\int_{\Omega_1}u^p(x,t) \varphi(x,t)\,dxdt\leq
    C R^{N+\alpha-\alpha(\ell-1)},
\end{equation}
where the constant $C$ on the right hand side of \rf{3.22} is independent
of $R$.  Note that
$N+\alpha-\alpha(\ell-1)\leq0$ if and only if  $ p\leq
1+{\alpha}/{N}.$ Now, we  consider two cases.

For $p<1+{\alpha}/{N},$ we have $N+\alpha-\alpha(\ell-1)<0.$
Hence, computing the
limit $R\to\infty$ in \rf{3.22} and using the
Lebesgue dominated convergence theorem, we obtain
$$
M_\infty=\int_{\mathbb{R}^N}u_0(x)\,dx
-\int_0^\infty\int_{\mathbb{R}^N}u^p(x,t)\,dxdt
\leq 2\ell\varepsilon\int_0^\infty\int_{\mathbb{R}^N}u^p\,dxdt.
$$
Since $u\in L^p(\mathbb{R}^N\times (0,\infty))$ (cf. \rf{3.2}) and
since
 $\varepsilon>0$ can be chosen arbitrary  small, we immediately
obtain that $M_\infty=0.$

In the critical case $p=1+{\alpha}/{N}$, we estimate first term
on the right hand side of inequality \rf{3.20} using again by  the
$\varepsilon$-Young inequality  and the second term by
 the H\"{o}lder inequality  (with $\bar p =p/(p-1)=\ell-1$)
as follows
\begin{eqnarray}
&&\hspace{-1cm}\int_{\Omega_1}u_0(x)\varphi(x,0)\,dx -\int_{\Omega_2}\int_{\Omega_1}u^p\varphi(x,t)\,dxdt\nonumber\\
&\leq&\ell\varepsilon\int_{\Omega_2}\int_{\Omega_1}u^p(x,t)\,dxdt \nonumber\\
&&+C(\varepsilon)\int_{\Omega_2}\int_{\Omega_1}\varphi_2^{\ell\bar{p}}(t)\varphi_1^{(\ell-1)\bar{p}}(x)\left|\Lambda^\alpha\varphi_1(x)\right|^{\bar{p}}\,dxdt\label{3.23}\\
&&+\ell\left(\int_{\Omega_3}\int_{\Omega_1}u^p(x,t)\,dxdt\right)^{1/p}\nonumber\\
&&\quad \times \left(\int_{\Omega_2}\int_{\Omega_1}\varphi_1^{\ell\bar{p}}(x)\varphi_2^{(\ell-1)\bar{p}}(t)\left|\partial_t\varphi_2(t)\right|^{\bar{p}}\,dxdt\right)^{1/\bar{p}}.\nonumber
\end{eqnarray}
Here,
$\Omega_3=\left\{t\in[0,\infty) \;:\; R^\alpha\leq t\leq
2R^\alpha\right\}$
 is the support of $\partial_t\varphi_2.$ Note that
$$
\int_{\Omega_3}\int_{\Omega_1}u^p(x,t)\,dx\,dt\to 0
\quad\hbox{as}\quad R\to\infty,$$
because $u\in L^p(\mathbb{R}^N\times[0,\infty))$ (cf. \rf{3.2}).

Now, introducing the new variables
$\xi=(BR)^{-1}x$, $\tau=R^{-\alpha}t$ and recalling that  $p=1+\alpha/N$,
we rewrite \rf{3.23} as follows
\begin{equation}\label{3.24}
\begin{split}
\int_{\Omega_1}u_0(x)\varphi(x,0)\,dx
&-\int_{\Omega_2}\int_{\Omega_1}u^p(x,t)\varphi(x,t)\,dxdt
-\varepsilon\ell\int_{\Omega_2}\int_{\Omega_1}u^p(x,t)\,dxdt\\
&\leq C_1B^{N/\bar p}\left(\int_{\Omega_3}\int_{\Omega_1}u^p(x,t)\,dxdt\right)^{1/p}+C_2C(\varepsilon)B^{-\alpha},
\end{split}
\end{equation}
where the constants $C_1, C2$ are independent of $R$, $B$, and of
$\varepsilon$. Passing in \rf{3.24} to the limit as $R\to+\infty$
and using the Lebesgue dominated convergence theorem we get
\begin{equation}\label{3.25}
\begin{split}
\int_{\mathbb{R}^N}u_0(x)\,dx&-\int_0^\infty\int_{\mathbb{R}^N}u^p(x,t)\,dxdt
-\varepsilon\ell \int_0^\infty\int_{\mathbb{R}^N}u^p(x,t)\,dxdt\\
&\leq C_2C(\varepsilon)B^{-\alpha}.
\end{split}
\end{equation}
Finally, computing the limit $B\to\infty$ in \rf{3.25}
we infer that $M_\infty=0$ beacuse
$\varepsilon>0$ can be arbitrarily small.
This complete the proof of Theorem 2. \qed

\bigskip

{\it Proof of Theorem 3.}
The proof proceeds by contradiction. Let $u$ be a non-negative
non-trivial  solution of \rf{eq}-\rf{ini}  with $\lambda=1$. Take
the test function $\varphi$ the same as in the proof of Theorem 2.
Repeating the estimations which lead to ($\ref{3.24}$), we obtain
\begin{equation}\label{3.26}
\begin{split}
&\int_{\Omega_1}u_0(x)\varphi(x,0)\,dx
+\int_{\Omega_2}\int_{\Omega_1}u^p(x,t) \varphi(x,t)\,dxdt\\
&\qquad\qquad-\varepsilon\ell\int_{\Omega_2}\int_{\Omega_1}u^p(x,t)\,dxdt\nonumber\\
&\leq C_1B^{N/\bar p}\left(\int_{\Omega_3}\int_{\Omega_1}u^p(x,t)\,dxdt\right)^{1/p}+ C_2C(\varepsilon)B^{-\alpha}.
\end{split}
\end{equation}
Now,  we chose $\varepsilon=1/(2\ell)$ in \rf{3.26} and we pass to
the following
limits: first $R\to\infty$, next $B\to\infty$.
Using the Lebesgue dominated convergence theorem, we obtain
$$\int_{\mathbb{R}^N}u_0(x)\,dx + \frac{1}{2}\int_0^\infty\int_{\mathbb{R}^N}u^p(x,t)\,dxdt\leq 0.
$$
Hence, $u(x,t)=0$ which contradicts our assumption imposed
on $u.$ \qed


\begin{thebibliography}{99}


\bibitem{BK99}
{\sc M.~Ben-Artzi \& H.~Koch}, {\it Decay of mass for a semilinear
  parabolic equation}, Comm. Partial Differential Equations
{\bf 24} (1999), 869--881.

\bibitem{BKW01} {\sc P. Biler, G. Karch \& W. A. Woyczy\'nski},
{\it Critical nonlinearity exponent and self-similar asymptotics
for L\'evy conservation laws},
Ann. I.H. Poincar\'e-Analyse non lin\'eare, {\bf 18} (2001), 613--637.

\bibitem{BLW02}
{\sc M. Birkner, J.A. Lopez-Mimbela, A.  Wakolbinger}, {\it Blow-up
   of semilinear PDE's at the critical dimension. A probabilistic approach.}
 Proc. Amer. Math. Soc. {\bf 130}
   (2002),  2431--2442.

\bibitem{5} {\sc A. C\'ordoba \& D. C\'ordoba},
{\it A maximum principle applied
to quasi-geostrophic equations}, Comm. Math. Phys. {\bf 249} (2004),
511--528.

\bibitem{6} {\sc J. Droniou \& C. Imbert},  {\it Fractal first-order
partial differential equations}, Arch. Rational Mech. Anal. {\bf 182}
(2006),  299--331.


\bibitem{8} {\sc H. Fujita}, {\it On the blowing up of solutions of the
problem for $u_t=\Delta u+u^{1+\alpha}$}, J. Fac. Sci. Univ. Tokyo
{\bf 13} (1966), 109--124.

\bibitem{9} {\sc M. Guedda \& M. Kirane},
{\it Criticality for some evolution
equations},   Differ. Uravn.  {\bf 37}  (2001),
511--520, 574-575.


\bibitem{10} {\sc M. Guedda \& M. Kirane},
{\it A note on nonexistence of global
solutions to a nonlinear integral equation},  Bull. Belg. Math. Soc.
Simon Stevin
{\bf 6} (1999),  491--497.

\bibitem{11} {\sc K. Hayakawa},  {\it On nonexistence of global solutions of
some semilinear parabolic differential equations},
Proc. Japan Acad. {\bf 49}
(1973), 503--505.


\bibitem{13} {\sc N. Ju},  {\it The maximum principle and the global attractor
for the dissipative 2D quasi-geostrophic equations}, Comm. Math.
Phys. {bf 255} (2005), 161--181.

\bibitem{KW08}
{\sc G. Karch \&  W.A. Woyczy\'nski},
{\it Fractal Hamilton-Jacobi-KPZ equations},
Trans. Amer. Math. Soc. {\bf 360} (2008), 2423--2442.

\bibitem{17} {\sc K. Kobayashi, T. Sirao \& H.  Tanaka},
{\it On the growing up problem for semilinear heat equations},  J. Math.
Soc. Japan {\bf 29} (1977),  407--424.

\bibitem{LS03}
{\sc Ph. Lauren\c cot \& Ph. Souplet}, {\it On the growth of mass for
a viscous Hamilton-Jacobi equation}, J. Anal.
Math. {\bf 89} (2003), 367--383.

\bibitem{19}  {\sc \`E. Mitidieri \& S. I. Pokhozhaev},
{\it Apriori estimates and the absence of solutions of nonlinear partial differential equations and
inequalities},
Tr. Mat. Inst. Steklova {\bf 234} (2001), 1--384; translation in Proc. Steklov Inst. Math. {\bf 234}
(2001), no. 3, 1--362.

\bibitem{20}  {\sc E. Mitidieri \& S. I. Pokhozhaev},  {\it Nonexistence
of weak solutions for some degenrate elliptic and parabolic problems
on ${\mathbb{R}}^N$}, J. Evol. Equ. {\bf 1} (2001), 189--220.

\bibitem{21}  {\sc R. G. Pinsky},  {\it Decay of mass for
the equation $u_t=\Delta u- a(x)u^p|\nabla u|^q$}, J. Diff.
Eq. {\bf 165} (2000), 1--23.

\bibitem{22}  {\sc S. Sugitani},  {\it On nonexistence of global solutions
for some nonlinear integral equations},  Osaka J. Math.  {\bf 12}
(1975), 45--51.

\end{thebibliography}
\end{document}